\documentclass[12pt,dvips]{amsart}

\usepackage[latin1]{inputenc}   
\usepackage[T1]{fontenc}        
\usepackage{amsmath,amsfonts,amssymb}
\usepackage{array}
\usepackage{rotating}
\usepackage{epsf}

\usepackage[matrix,frame]{xypic}

\newtheorem{example}{Example}[section]

\newtheorem{theorem}[example]{Theorem}

\newtheorem{corollary}[example]{Corollary}

\newtheorem{definition}[example]{Definition}
\newtheorem{proposition}[example]{Proposition}

\newtheorem{lemma}[example]{Lemma}

\def\Proof{\noindent \it Proof -- \rm}
\def\qed{\hspace{3.5mm} \hfill \vbox{\hrule height 3pt depth 2 pt width 2mm}
\bigskip}

\def\NCSF{{\bf Sym}}           

\def\inv{{\rm inv}}     

\def\ssh{\Cup}          

\def\<{\langle}
\def\>{\rangle}
\def\NN{{\mathbb N}}    
\def\RR{{\mathbb R}}    
\def\KK{{\mathbb K}\, } 

\def\SG{{\mathfrak S}}  

\def\Des{\operatorname{Des}}

\def\shuff#1#2{\mathbin{
\hbox{\vbox{ \hbox{\vrule \hskip#2 \vrule height#1 width 0pt
}%
\hrule}%
\vbox{ \hbox{\vrule \hskip#2 \vrule height#1 width 0pt
\vrule }%
\hrule}%
}}}

\def\shuf{{\mathchoice{\shuff{7pt}{3.5pt}}%
{\shuff{6pt}{3pt}}%
{\shuff{4pt}{2pt}}%
{\shuff{3pt}{1.5pt}}}}%
\def\shuffle{\,\shuf\,}

\def\id{\operatorname{id}}
\def\A{{\bf A}}
\def\des{\operatorname{des}}
\def\maj{\operatorname{maj}}
\def\inv{\operatorname{inv}}

\def\vect{{}}
\def\Invcode{\operatorname{Ic}}
\def\Majcode{\operatorname{Mc}}
\def\Lcode{\operatorname{Lc}}
\def\Scode{\operatorname{Sc}}
\def\Gcode{\operatorname{Gc}}
\def\Card{\operatorname{Card}}
\def\eval{\operatorname{ev}}
\def\X{{\bf X}}

\title[Generalizations of the Foata-Sch\"utzenberger equidistribution]%
{Multivariate generalizations \\
       of the Foata-Sch\"utzenberger \\
       equidistribution}

\author[F. Hivert, J.-C.~Novelli, and J.-Y.~Thibon]%
{Florent Hivert, Jean-Christophe Novelli, and Jean-Yves Thibon}

\address[Hivert]{LITIS, Universit\'e de Rouen ; Avenue de l'universit\'e ;
76801 Saint \'Etienne du Rouvray, France\\}

\address[Novelli and Thibon]{Institut Gaspard Monge, Universit\'e de
Marne-la-Vall\'ee \\
5, Boulevard Descartes \\Champs-sur-Marne \\77454 Marne-la-Vall\'ee cedex 2 \\
FRANCE}
\email[Florent Hivert]{hivert@univ-mlv.fr}
\email[Jean-Christophe Novelli]{novelli@univ-mlv.fr}
\email[Jean-Yves Thibon]{jyt@univ-mlv.fr} 

\begin{document}

\begin{abstract}
A result of Foata and Sch\"utzenberger states that two statistics on
permutations, the number of inversions and the inverse major index, have the
same distribution on a descent class. We give a multivariate generalization of
this property: the sorted vectors of the Lehmer code, of the inverse majcode,
and of a new code (the inverse saillance code), have the same
distribution on a descent class, and their common multivariate generating
function is a flagged ribbon Schur function.
\end{abstract}

\maketitle

\section{Introduction}

The \emph{major index} of a permutation, discovered by Major Percy
Alexander MacMahon, and named after his military rank, is the sum of its
\emph{descents}
\begin{equation}
\maj(\sigma) = \sum_{\sigma(i)>\sigma(i+1)} i.
\end{equation}
The maximum value of $\maj$ over the set $\SG_n$ of permutations of size $n$
is $n(n-1)/2$, the same as for the inversion number $\inv$, and MacMahon
proved~\cite{MM} (actually, he proved a similar result for an arbitrary
rearrangement class of words, but in this paper, we will only deal with
permutations) that both statistics have the same distribution
\begin{equation}
\sum_{\sigma\in\SG_n} q^{\maj(\sigma)} =
\sum_{\sigma\in\SG_n} q^{\inv(\sigma)} = [n]_q! =
\prod_{i=1}^n \frac{1-q^i}{1-q}.
\end{equation}
More than fifty years later, Foata and Sch\"utzenberger~\cite{FS} proved
that equidistribution holds on a \emph{descent class}:
\begin{equation}
\label{FS}
\sum_{\Des(\sigma)=D} q^{\maj(\sigma^{-1})} =
\sum_{\Des(\sigma)=D} q^{\inv(\sigma)}
\end{equation}
where $\Des(\sigma) = \{ i\,|\, \sigma(i)>\sigma(i+1) \}$ is the descent set
of $\sigma$.
The original proof (and, up to recently the only one, \emph{cf.}~\cite{Loth2},
chapter~11) of this result was bijective.
In this note, we obtain a multivariate refinement of~(\ref{FS}) (different
from the one of~\cite{FH}): we prove
that, up to order, the three integer vectors to be defined below, namely
the inverse Lehmer code $\Invcode(\sigma^{-1})=(c_1,\ldots,c_n)$, the inverse
major code $\Majcode(\sigma^{-1})=(m_1,\ldots,m_n)$, and a new code
$\Scode(\sigma^{-1})=(s_1,\ldots,s_n)$ have the same distribution on a
descent class, that is, if $x_0,\ldots,x_{n-1}$ are independent indeterminates
\begin{equation}
\label{Nous}
\sum_{\Des(\sigma)=D} \prod_i x_{c_i(\sigma^{-1})} =
\sum_{\Des(\sigma)=D} \prod_i x_{m_i(\sigma^{-1})} =
\sum_{\Des(\sigma)=D} \prod_i x_{s_i(\sigma^{-1})}.
\end{equation}
Indeed, the codes are defined in such a way that
\begin{equation}
\sum_{i=1}^n c_i(\sigma) = \inv(\sigma),\quad
\sum_{i=1}^n m_i(\sigma^{-1}) = \maj(\sigma^{-1}),
\end{equation}
so that Equation~(\ref{Nous}) is a refinement of Equation~(\ref{FS}).
\goodbreak

{\footnotesize
{\it Acknowledgements.-}
This project has been partially supported by CNRS.
The authors would also like to thank the contributors of the MuPAD project,
and especially those of the combinat package, for providing the development
environment for this research (see~\cite{HT} for an introduction to
MuPAD-Combinat).
}

\section{Notations}

\subsection*{Alphabets and operations on words}

In all the paper, we deal with a totally ordered infinite alphabet $A$,
represented either by $\{a,b,c,\ldots\}$ or by $\{1,2,3,\ldots\}$.
The \emph{free associative algebra} over $A$
is denoted by $\KK\<A\>$, where $\KK$ is some field of characteristic zero.
The \emph{evaluation} $\eval(w)$ of a word $w$ of size $n$ over the alphabet
$\{0,\ldots,n\}$ is the list of numbers of appearance $\Card\{i|w_i=a\}$
of all letters $a\in A$ in $w$. For example, the evaluation of $45143251812$
is $032122001000$.
We denote by $\id_n=12\cdots n$ the identity permutation of size $n$.

\noindent
The \emph{shuffle product} $w_1\shuffle w_2$  of two words $w_1$ and $w_2$
is recursively defined by
$w_1 \shuffle \epsilon  = w_1$ and $\epsilon \shuffle w_2  = w_2$, where
$\epsilon$ is the empty word, and
\begin{equation}
au \shuffle bv = a(u \shuffle bv) + b(au \shuffle v),
\qquad \text{$a$, $b\in A$,\quad $u$, $v\in A^*$.}
\end{equation}
For example,
\begin{equation}
12\shuffle 43 = 1243 + 1423 + 1432 + 4123 + 4132 + 4312\,.
\end{equation}
For a word $w=w_1\dots w_n$ over the integers, and $k\in\NN$, we denote by
$w[k]$ the \emph{shifted word}
\begin{equation}
w[k] := (w_1+k)\cdot(w_2+k)\cdots (w_n+k).
\end{equation}
The \emph{shifted shuffle} of two permutations $\alpha\in\SG_k$ and
$\beta\in\SG_l$ is then defined by
\begin{equation}
\alpha \ssh \beta := \alpha \shuffle (\beta[k]).
\end{equation}
%
%

\subsection*{Compositions}

A \emph{composition} of an integer $n$ is a sequence of positive integers
of sum $n$.
The \emph{descent set} $\Des(I)$ of a composition $I=(i_1,\ldots,i_r)$ is the
set of partial sums $\{i_1,i_1+i_2,\ldots,i_1+\cdots +i_r\}$.
Compositions are ordered by $I\leq J$ iff $\Des(I)\subseteq\Des(J)$. In this
case, we say that $I$ is \emph{coarser} than $J$.

The \emph{descent composition} $I=C(\sigma)$ of a permutation $\sigma$ is the
composition of $n$ whose descents are equal to the descents of $\sigma$, that
is, the set of integers $j$ such that $\sigma(j)>\sigma(j+1)$.

If $I=(i_1,\ldots,i_r)$ is a composition of $n$, let $D_{\leq I}$ be the sum
of all permutations whose descent composition is
coarser than $I$.
Then
\begin{equation}
D_{\leq I} = (\id_{i_1}\ssh\id_{i_2} \ssh \cdots \ssh \id_{i_r})^\vee
\end{equation}
where ${}^\vee$ is the linear involution sending each permutation to its
inverse.
The sum of all permutations whose descent composition is $I$ will be denoted
by $D_I$.
\smallskip

Recall that the algebra of noncommutative symmetric functions $\NCSF$ is the
free associative algebra on symbols $S_n$ so that a basis is given by the
$S^I=S_{i_1}\cdots S_{i_r}$ for all compositions
$I=(i_1,\ldots,i_r)$~\cite{NCSF1}. We will make use of this basis of $\NCSF$
and of the ribbon basis $R_I$ defined by
\begin{equation}
S^I := \sum_{J\leq I} R_J.
\end{equation}
When $A$ is an ordered alphabet, $S_n(A)$ can be realized as the sum of all
nondecreasing words in $A^n$.
The commutative image of $\NCSF$ is the algebra of symmetric functions.
The $S_n$ are mapped to the usual complete homogeneous functions $h_n$, and
the $R_I$ to the ribbon Schur functions $r_I$.

\subsection*{Codes}

Let us say that a sequence $a=(a_1,\ldots,a_n)$ is
\emph{sub-diagonal} if $0\leq a_i\leq n-i$. 
A \emph{code} is then a bijection between the symmetric group $\SG_n$ and the
set of sub-diagonal sequences of length $n$. Among known codes, we will be
interested in the \emph{Lehmer code} and the \emph{major code}.

Recall that the \emph{Lehmer code} (or \emph{Lcode}, for short)
$\Lcode(\sigma)$ of a permutation $\sigma\in\SG_n$ is the sequence
$(c_i)_{1\leq i\leq n}$, where
\begin{equation}
c_i = \Card\{j>i\ |\ \sigma_j<\sigma_i\}.
\end{equation}
For example, the code of the permutation $531962487$ is $420520010$.

The Lehmer code of the inverse permutation will be called the
\emph{Inversion code} (or \emph{invcode}, for short).
It is the sequence $\Invcode(\sigma)=(a_1,\ldots,a_n)$, where $a_i$ is
the number of values greater than $i$ to its left.
In other words, the invcode, as the Lehmer code, splits the inversions of
$\sigma$ into blocks.

Since the number of inversions of a permutation is the sum of the components
of its code, one may look for an analogous vector having as sum the major
index. The \emph{major code} (or \emph{majcode}, for short) solves this
question. It is implicit in Carlitz~\cite{Car} and explicitly stated by
Rawlings in~\cite{Raw} (see also~\cite{Ska}).
Recall that the \emph{major index} $\maj$ of a permutation is the sum of the
positions of its descents.
Now, for $\sigma\in\SG_n$, denote by $\sigma^{(i)}$ the subword of $\sigma$
obtained by erasing the letters smaller than $i$, so that
$\sigma=\sigma^{(1)}$.
Then the majcode $\Majcode(\sigma)$ of $\sigma$ is the sequence
$(c_i)_{1\leq i\leq n}$, where $c_n=0$ and
\begin{equation}
c_i = \maj(\sigma^{(i)}) - \maj(\sigma^{(i+1)}),
\end{equation}
for all $1\leq i\leq n-1$.
For example, $\Majcode(935721468)=501012010$.

Finally, the \emph{sorted vector} $V^\uparrow$ of a vector
$(v_1,\ldots,v_n)$ is its nondecreasing rearrangement.

\section{Cayley trees and codes}

\subsection{From differential equations to trees}

Cayley~\cite{Cay} introduced trees in order to solve the differential equation
\begin{equation}
\label{vit}
\frac{dx}{dt}(t) = {\vect V}(x(t)),
\end{equation}
where $\vect V$ is a vector field, that is a function from $\RR^d$ to itself.
 
Formally, the special case $d=1$ gives the following values for the
coefficients $x_n=d^nx/dt^n(0)$ of the Taylor expansion at $t=0$ of the
solution:
\begin{alignat}{1}
x_1 =\ & V_0\\ 
x_2 =\ & V_1 V_0\\ 
x_3 =\ & V_2 V_0^2 + V_1^2 V_0\\ 
x_4 =\ & V_3 V_0^3 + 4\, V_2 V_1 V_0^2 + V_1^3 V_0\\ 
x_5 =\ & V_4 V_0^4 + 7\, V_3 V_1 V_0^3 + 4\, V_2^2 V_0^3 + 11\,
V_2 V_1^2 V_0^2 + V_1^4 V_0\\ 
x_6 =\ & V_5 V_0^5 + 11\, V_4 V_1 V_0^4 + 15\, V_3 V_2 V_0^4 +
32\, V_3 V_1^2 V_0^3 + 34\, V_2^2 V_1 V_0^3 \\
\notag & \hspace{6cm} + 26\, V_2 V_1^3 V_0^2 + V_1^5 V_0
\end{alignat}
where $V_n$ is $\frac{d^n V}{dx^n}(x(0))$.

Assuming without loss of generality that $V_0=1$, we have in the
one-dimensional case
\begin{equation}
x_{n+1} = C_n(V_1,\ldots,V_n),
\end{equation}
where the polynomials $C_n$ reduce to the Eulerian polynomials
\begin{equation}
C_n(q,\ldots,q) = A_n(q)
\end{equation}
when all the $V_i$ are equal to $q$.

The polynomials $C_n$ giving the Taylor coefficients of the unique solution of
$\frac{dx}{dt} = V(x(t))$ with $x(0)=0$ and $V(0)=1$ should be compared to the
exponential Bell polynomials
\begin{equation}
  B_n(x_1,\ldots,x_n) = \sum_k B_{n,k}(x_1,\dots x_{n+1-k})
\end{equation}
giving the Taylor coefficients of $y(t)=V(x(t))$
\begin{equation}
y_n = \sum_{k=0}^{n} V_k B_{n,k}(x_1,\dots,x_{n+1-k}).
\end{equation}
Both calculations are related to the so-called Faa di Bruno Hopf algebra
\cite{CHAlg}. The $C_n$ are to the Eulerian polynomials what the $B_n$ are to
the (one variable) Bell polynomials
\begin{equation}
  b_n(q) = B_n(q, \dots, q)\,.
\end{equation}

It is immediate that the sum of the coefficients of the monomials $V^\alpha$
in $x_{n+1}$ in $n!$, so that the exponent vector $\alpha$ should be
interpretable as a multivariate statistic on the symmetric group $\SG_n$.
Alain Lascoux observed that this statistic seems to coincide with the sorted
evaluation of the Lehmer code, which is fundamental in the theory of Schubert
and Grothendieck polynomials~\cite{Las}.
The codes of the permutations of $\SG_3$ are
\begin{equation}
000,\ 010,\ 100,\ 110,\ 200,\ 210,
\end{equation}
whose evaluations are
\begin{equation}
3000,\ 2100,\ 2100,\ 1200,\ 2010,\ 1110,
\end{equation}
giving back the coefficients of $x_4$.
More generally, one has
\begin{equation}
x_n = \sum_{\sigma\in\SG_{n-1}} V^{\Invcode(\sigma)}.
\end{equation}

We shall see that a better way to understand this formula is to rely upon
another statistic on permutations coming up more naturally than the Lehmer
code from differential equations, namely the \emph{saillance code}.
It comes from decreasing tree structures on permutations
appearing in the $d$-dimensional case of Equation~(\ref{vit}).
When $d>1$, the derivatives $V_i$ have to be replaced by the differentials
$D^k V$ defined by
\begin{equation}
[D^k \vect V (\vect U_1, \dots \vect U_k)]_i:=
\sum_{j_1 \dots j_k = 1}^d
\frac{\partial^k [\vect V]_i}{\partial x_{j_1} \dots \partial x_{j_k}}
[\vect U_1]_{j_1} \dots [\vect U_k]_{j_k}\,,
\end{equation}
where $\vect U_j$ is any vector field and $[\vect U_j]_i$ denotes the $i$-th
coordinate of $\vect U_j$.
If the $U_i$ are evaluated at $x(t)$,
\begin{equation}
\frac{d (D^k \vect V (\vect U_1, \dots, \vect U_k))}{dt}\ =\
D^{k+1} \vect V (\vect U_1, \dots, \vect U_k, V)\ +\
\sum_{i=1}^{k} D^k \vect V (\vect U_1, \dots, \frac{d U_i}{dt},\dots, \vect
U_k)\,.
\end{equation}
The calculation of the Taylor coefficients $x_{n+1}=d^n/dt^n V(x(t))$ at $t=0$
gives rise to nested derivatives, which can be conveniently represented by
topological trees.

Define $D_T \vect V$ recursively by:
if the root of $T$ has $k$ subtrees $T_1,\dots,T_k$,
\begin{equation}
D_T \vect V := D^k \vect V (D_{T_1} \vect V, \dots, D_{T_k} \vect V)\,,
\end{equation}
with the convention $D^0 V = V$.

For example,
\begin{equation}
  \entrymodifiers={+<2pt>[F]}
  D^2 \vect V(\vect V, D^3 \vect V(\vect V, D^2 \vect V(\vect V, \vect V),
\vect V))\ =\
  \vcenter{\xymatrix@C=2mm@R=2mm{
   *{}& {D^2 \vect V}\ar@{-}[drr]\ar@{-}[dl] \\
   {\vect V} &*{}&*{}& {D^3 \vect V}\ar@{-}[drr]\ar@{-}[d]\ar@{-}[dll] \\
   *{}& {\vect V} &*{}& {D^2 \vect V}\ar@{-}[dr]\ar@{-}[dl] &*{}& {\vect V} \\
   *{}&*{}& {\vect V} &       *{}           & {\vect V}
  }}
\ =\
  \entrymodifiers={+<4pt>[o][F-]}
  D_{\vcenter{\xymatrix@C=2mm@R=2mm{
   *{}& *+<4pt>[o][F**]{}\ar@{-}[drr]\ar@{-}[dl] \\
   {} &*{}&*{}& {}\ar@{-}[drr]\ar@{-}[d]\ar@{-}[dll] \\
   *{}& {} &*{}& {}\ar@{-}[dr]\ar@{-}[dl] &*{}& {} \\
   *{}&*{}& {} &       *{}           & {}
  }}}V.
\end{equation}
With this notation, the derivation rule becomes
\begin{equation}
\frac{d (D_T V)}{dt} = \sum_{T'} D_{T'} V,
\end{equation}
where $T'$ runs over the set of trees obtained from $T$ by adding a leaf
to each node of $T$.
For example, writing only the trees in the previous equation, one has
\entrymodifiers={+<4pt>[o][F-]}
\begin{equation}
\frac{d(\vcenter{\xymatrix@C=1mm@R=1mm{
        *{} & *+<4pt>[o][F**]{}\ar@{-}[dl]\ar@{-}[dr] \\ 
            & *{} & {}  
      }})}{dt}
\ =\ 
     \vcenter{\xymatrix@C=1mm@R=1mm{
        *{} & *+<4pt>[o][F**]{}\ar@{-}[dl]\ar@{-}[d]\ar@{-}[dr] \\ 
            &    & {}
      }}
    \ +\ 
    \vcenter{\xymatrix@C=1mm@R=1mm{
        *{} & *+<4pt>[o][F**]{}\ar@{-}[dl]\ar@{-}[dr] \\ 
         \ar@{-}[d] & *{} &  \\
        {}  
      }}
    \ + \ 
    \vcenter{\xymatrix@C=1mm@R=1mm{
        *{} & *+<4pt>[o][F**]{}\ar@{-}[dl]\ar@{-}[dr] \\ 
            & *{} & \ar@{-}[d] \\
        *{} & *{} &  {}  
      }}
\ =\ 
     \vcenter{\xymatrix@C=1mm@R=1mm{
        *{} & *+<4pt>[o][F**]{}\ar@{-}[dl]\ar@{-}[d]\ar@{-}[dr] \\ 
            &    & {}
      }}
    \ +\ 2\, 
    \vcenter{\xymatrix@C=1mm@R=1mm{
        *{} & *+<4pt>[o][F**]{}\ar@{-}[dl]\ar@{-}[dr] \\ 
            & *{} & \ar@{-}[d] \\
        *{} & *{} &  {}  
      }}
\end{equation}
We have thus, for each integer $n$, an expression
\begin{equation}
\label{bla}
x_{n} = \sum_{|T|=n} \alpha_T D_TV,
\end{equation}
where the $\alpha_T$ are positive integers, sometimes known as the
Connes-Moscovici coefficients~\cite{CM,Chap}.

The first values are
\entrymodifiers={+<4pt>[o][F-]}
\begin{alignat}{1}
x_1\    =\ \ & \bullet \\
x_2\    =\ \ & 
    \vcenter{\xymatrix@C=1mm@R=1mm{
        *+<4pt>[o][F**]{}\ar@{-}[d] \\ {}  
      }} 
    \\
x_3\    =\ \ &
    \vcenter{\xymatrix@C=1mm@R=1mm{
        *+<4pt>[o][F**]{}\ar@{-}[d] \\ \ar@{-}[d] \\ {}  
      }} 
    +
    \vcenter{\xymatrix@C=1mm@R=1mm{
        *{} & *+<4pt>[o][F**]{}\ar@{-}[dl]\ar@{-}[dr] \\ 
            & *{} & {}  
      }} 
    \\
x_4\    =\ \ &
    \vcenter{\xymatrix@C=1mm@R=1mm{
        *{} & *+<4pt>[o][F**]{}\ar@{-}[dl]\ar@{-}[d]\ar@{-}[dr] \\ 
            &    & {}
      }}
    + 3 \ 
    \vcenter{\xymatrix@C=1mm@R=1mm{
        *{} & *+<4pt>[o][F**]{}\ar@{-}[dl]\ar@{-}[dr] \\ 
            & *{} & \ar@{-}[d] \\
        *{} & *{} &  {}  
      }}
    +
    \vcenter{\xymatrix@C=1mm@R=1mm{
        *{} & *+<4pt>[o][F**]{}\ar@{-}[d] \\ 
        *{} & \ar@{-}[dl]\ar@{-}[dr]  \\
           & *{} & {}   
      }}
    +
    \vcenter{\xymatrix@C=1mm@R=1mm{
        *+<4pt>[o][F**]{}\ar@{-}[d] \\ \ar@{-}[d] \\ \ar@{-}[d] \\ {}
      }} 
    \\
x_5\    =\ \ &
    \vcenter{\xymatrix@C=0.5mm@R=2mm{
        *{}&*{}&*{}&*=<0pt>{\bullet}\ar@{-}[dlll]\ar@{-}[dl]\ar@{-}[dr]\ar@{-}[drrr]\\ 
           &*{}&   &*{}&   &*{}& {}
      }}
    + 6 \ 
    \vcenter{\xymatrix@C=1mm@R=1mm{
        *{} & *+<4pt>[o][F**]{}\ar@{-}[dl]\ar@{-}[d]\ar@{-}[dr] \\ 
            &     & \ar@{-}[d] \\
        *{} & *{} &  {}  
      }}
    + 4 \ 
    \vcenter{\xymatrix@C=1mm@R=1mm{
        *{} & *+<4pt>[o][F**]{}\ar@{-}[dl]\ar@{-}[dr] \\ 
            & *{} &\ar@{-}[dl]\ar@{-}[dr] \\
        *{} &     & *{} & {}  
      }}
    + 4 \ 
    \vcenter{\xymatrix@C=1mm@R=1mm{
        *{} & *+<4pt>[o][F**]{}\ar@{-}[dl]\ar@{-}[dr] \\ 
            & *{} &\ar@{-}[d] \\
        *{} & *{} &\ar@{-}[d] \\
        *{} & *{} & {}  
      }}
    + 3 \ 
    \vcenter{\xymatrix@C=1mm@R=1mm{
        *{} & *+<4pt>[o][F**]{}\ar@{-}[dl]\ar@{-}[dr] \\ 
        \ar@{-}[d] & *{} &\ar@{-}[d] \\
        {} & *{} & {}  
      }}
+    \vcenter{\xymatrix@C=1mm@R=1mm{
        *{} & *+<4pt>[o][F**]{}\ar@{-}[d] \\ 
        *{} & \ar@{-}[dl]\ar@{-}[d]\ar@{-}[dr] \\ 
            &    & {}
      }}
    + 3 \ 
    \vcenter{\xymatrix@C=1mm@R=1mm{
        *{} & *+<4pt>[o][F**]{}\ar@{-}[d] \\ 
        *{} & \ar@{-}[dl]\ar@{-}[dr] \\ 
            & *{} & \ar@{-}[d] \\
        *{} & *{} &  {}  
      }}
    +
    \vcenter{\xymatrix@C=1mm@R=1mm{
        *{} & *+<4pt>[o][F**]{}\ar@{-}[d] \\
        *{} & \ar@{-}[d] \\ 
        *{} & \ar@{-}[dl]\ar@{-}[dr]  \\
           & *{} & {}   
      }}
    +
    \vcenter{\xymatrix@C=1mm@R=1mm{
        *+<4pt>[o][F**]{}\ar@{-}[d] \\
        \ar@{-}[d] \\ \ar@{-}[d] \\ \ar@{-}[d] \\ {}
      }} 
\end{alignat}

\subsection{From trees to permutations and statistics}

It is easy to see that the coefficient $\alpha_T$ is equal to the number of
\emph{increasing trees} of shape $T$, that is, the set of trees obtained by
labelling the nodes of $T$ by the integers from $1$ to $n$ so that the
label of each node is greater than the label of its father.
Recall that increasing topological trees have a \emph{canonical form} that
consists in ordering the children of each node so that they are increasing
from left to right.
For example, written in canonical form, here are the five increasing trees of
that given shape.
\entrymodifiers={+<4pt>}
\begin{equation}
\label{canonik}
    \vcenter{\xymatrix@C=2mm@R=2mm{
        *{}  & {1}\ar@{-}[dr]\ar@{-}[dl] \\
        {2}  & *{} & {3}\ar@{-}[d] \\ 
        *{}  & *{} & {4}\ar@{-}[dl]\ar@{-}[dr]  \\
        *{}  & {5} & *{} & {6}   
      }}
    \qquad
    \vcenter{\xymatrix@C=2mm@R=2mm{
        *{}  & *{}  & {1}\ar@{-}[dr]\ar@{-}[dl] \\
        *{}  & {2}\ar@{-}[d]  & *{} & {3} \\ 
        *{} & {4}\ar@{-}[dl]\ar@{-}[dr]  \\
        {5} & *{} & {6}   
      }}
    \qquad
    \vcenter{\xymatrix@C=2mm@R=2mm{
        *{}  & *{}  & {1}\ar@{-}[dr]\ar@{-}[dl] \\
        *{}  & {2}\ar@{-}[d]  & *{} & {4} \\ 
        *{} & {3}\ar@{-}[dl]\ar@{-}[dr]  \\
        {5} & *{} & {6}   
      }}
    \qquad
    \vcenter{\xymatrix@C=2mm@R=2mm{
        *{}  & *{}  & {1}\ar@{-}[dr]\ar@{-}[dl] \\
        *{}  & {2}\ar@{-}[d]  & *{} & {5} \\ 
        *{} & {3}\ar@{-}[dl]\ar@{-}[dr]  \\
        {4} & *{} & {6}   
      }}
    \qquad
    \vcenter{\xymatrix@C=2mm@R=2mm{
        *{}  & *{}  & {1}\ar@{-}[dr]\ar@{-}[dl] \\
        *{}  & {2}\ar@{-}[d]  & *{} & {6} \\ 
        *{} & {3}\ar@{-}[dl]\ar@{-}[dr]  \\
        {4} & *{} & {5}   
      }}
\end{equation}

Recall that the total number of increasing trees of size $n$ is $(n-1)!$ There
exist many different bijections with permutations of size $n-1$.
We shall make use of the following one.
Start from an increasing tree of size $n$, replace each label $l$ in $T$ by
$n+1-l$, then reorder the (now decreasing) tree in canonical form and read it
in prefix order, forgetting the root.
For example, the permutations corresponding with the increasing trees
of~(\ref{canonik}) are

\begin{equation}
\label{perms}
43125,\qquad\qquad  45312,\qquad\qquad  35412,\qquad\qquad  25413,\qquad\qquad  15423.
\end{equation}

The special case $d=1$ amounts to replace each tree $T$ by the monomial
\begin{equation}
\prod_{o \in \operatorname{nodes}(T)} V_{\operatorname{arity}(o)}.
\end{equation}
This statistic is the evaluation of a code.

\begin{definition}
The \emph{saillance code} (or \emph{scode}, for short) $\Scode(T)$ of a tree
$T$ of size $n$ as the sequence of labels of the fathers (minus one) of
$n, n-1, \ldots, 2$.
\end{definition}

For example, the scode of the trees of~(\ref{canonik}) are
\begin{equation}
\label{scodes}
33200,\qquad \qquad 33100,\qquad \qquad 22010,\qquad \qquad 20210,\qquad
\qquad 02210.
\end{equation}

\section{Properties of the scode}

The scode can be directly defined on permutations as follows:
the scode of a permutation $\sigma\in\SG_n$ is the sequence
$a=(a_1,\ldots,a_n)$, where $a_i$ is the number of letters of $\sigma$ greater
than or equal to the rightmost letter to the left of $i$ and greater than
$i$.

For example, one can check that the scodes of the permutations
of~(\ref{perms}) are the sequences of~(\ref{scodes}).

\begin{proposition}
The scode is a code.
\end{proposition}

\Proof
Since $a_i$ counts a number of letters in $\sigma$ all greater than $i$,
the sequence $a$ is sub-diagonal. Moreover, thanks to the bijection between
permutations and increasing trees, and to the fact that the scode is obviously
injective from trees, the scode is a bijection.
\qed

The algorithm giving back the permutation from its scode
$a=(a_1,\ldots,a_n)$ is as follows: put $n$ and then insert letters
$n-1$ to $1$ such that $i$ is inserted immediately after letter $n+1-a_i$
(and first if $a_i=0$).

\begin{proposition}
Interpreting a sub-diagonal sequence as the scode of a
permutation, the number of descents of this permutation is given by the number
of non-zero different values in its scode.
\end{proposition}

\Proof
The number of descents of a permutation is equal to the number of internal
nodes except the root of its corresponding increasing tree. This
last number is obviously equal to the number of non-zero different values of
its scode.
\qed

It is known that the generating function of the Lehmer code (up to order) over
a class $D_{\leq I}$ admits a closed expression, as a product of complete
symmetric functions over a flag of alphabets.
This property is also true for the scode.
Let us first introduce for all $n\geq0$, the alphabet
$X_n= \{x_0,\ldots,x_n\}$ where the $x_i$ are commuting indeterminates.
With a given sub-diagonal sequence $c$, we associate the monomial
\begin{equation}
x_c = x_{c_1} x_{c_2} \cdots x_{c_n} \in X_{n-1}\times X_{n-2}
\times\cdots\times X_0.
\end{equation}

\begin{theorem}
\label{scode-thsg}
Let $I=(i_1,\ldots,i_r)$ be a composition.
The sum of the scodes of the inverses of the elements of $D_{\leq I}$ are
given by the generating function
\begin{equation}
\begin{split}
\label{sg-Sc}
F_S(I) := \sum_{\sigma\in D_{\leq I}} x_{\Scode(\sigma^{-1})} &=
\sum_{\sigma\ \in\ \id_{i_1}\!\!\ssh\cdots\ssh\, \id_{i_r}} x_{\Scode(\sigma)}
\\ &= h_{i_1}(X_{i_2+\cdots+i_r})\ 
h_{i_2}(X_{i_3+\cdots+i_r}) \cdots h_{i_{r-1}}(X_{i_r})\  h_{i_r}(X_0).
\end{split}
\end{equation}
The right-hand side will be denoted by $h^I(\X_I)$, where $\X_I$ denotes the
\emph{flag of alphabets}
\begin{equation}
(X_{n-i_1},X_{n-i_1-i_2},\ldots,X_{n-n}).
\end{equation}
\end{theorem}

The proof relies on the following lemmas.

\medskip
Let $\beta\in\SG_n$. For $0\leq i\leq n$, denote by $1\ssh_i \beta$ the
term of the shifted shuffle $1\ssh \beta$ in which $1$ occurs at the
$(i+1)$-st position, \emph{e.g.},
$1\ssh_0 21 = 132;\quad
1\ssh_1 21 = 312;\quad
1\ssh_2 21 = 321.$

\begin{lemma}
Let $\beta\in\SG_n$. Then
\begin{equation}
x_{\Scode(1\ssh_i \beta)} = x_{\tau_S(\beta)(i)} x_{\Scode(\beta)},
\end{equation}
where $\tau_S(\beta)$ is the permutation of $\{0,\ldots,n\}$ defined by
\begin{equation}
\label{tau-sco}
\begin{split}
\tau_S(\beta)(0)=0,
\qquad\text{and}\qquad
\tau_S(\beta)(i)=n+1-\beta(i).
\end{split}
\end{equation}
\end{lemma}

For example,
\begin{equation}
\begin{split}
\beta &=\,\, 9\ 4\ 1\ 6\ 2\ 5\ 7\ 3\ 8 \\
\tau_S(\beta) &= 0\ 1\ 6\ 9\ 4\ 8\ 5\ 3\ 7\ 2
\end{split}
\end{equation}

\Proof
Let $\beta'=1\ssh_i\beta$. It is obvious that
$\Scode_{i+1}(\beta')=\Scode_i(\beta)$ for $i\in [1,n]$.
So $x_{\Scode(\beta')}/x_{\Scode(\beta)} = x_{\Scode_1(\beta')}$.
The value of the rightmost value to the left of $1$ and greater than $1$ is
its neighbour to the left. And the number of values greater than or equal to
this last value is its complement to $n+1$, so that it corresponds to the
definition of $\tau$.
\qed

%
\begin{lemma}
\label{sco-l3}
Let $\beta\in\SG_n$, let $k$ be an element of $[0,n]$
and let $\beta'=1\ssh_k\beta$.
Then
\begin{equation}
\label{egalite}
\tau_S({\beta'})(i) = \tau_S({\beta})(i)
\end{equation}
for $i\in[0,k]$.
\qed
\end{lemma}

For example, given $\beta=72451836$, we have $\tau_S(\beta)= 027548163$.
The case $k=3$ gives $\beta'=835162947$ and $\tau_S({\beta'})=0275948163$.
The case $k=7$ gives $\beta'=835629417$ and $\tau_S({\beta'})=0275481693$.

\begin{lemma}
Let $\beta\in\SG_n$. Then the scodes of the elements in $id_k\ssh\beta$ are
\begin{equation}
(\tau_S(\beta)(i_1),\tau_S(\beta)(i_2),\ldots,
 \tau_S(\beta)(i_k),\Scode(\beta))
\end{equation}
where $(i_j)$ runs over nondecreasing sequences in $[0,n]$,
$0\leq i_1\leq i_2\leq \ldots\leq i_k\leq n$.

In particular, we have
\begin{equation}
\label{scode-X}
\sum_{\sigma\in\id_k \ssh\beta} x_{\Scode(\sigma)}
= h_k(X_n)\ x_{\Scode(\beta)}.
\end{equation}
\end{lemma}

\Proof
It is sufficient to prove the result for $k=2$. The computation of
$12\ssh\beta$ can be decomposed as the shifted shuffle of $1$ with $\beta$
followed by the shifted shuffle of $1$ with the new elements where $1$ cannot
go to the right of $2$. Since the set of values $\tau_S({\beta'})(i)$ for all
$i$ in $[0,k]$ where $k$ is the position of $1$ in $\beta'$ is equal to the
set of values $\tau_S({\beta})(i)$ for all $i$ in $[0,k]$, we are done for the
first part of the lemma.

Then, since $\tau$ is a permutation of $[0,n]$, the commutative
image of all words
$a_{\tau_\beta(i_1)}\ldots a_{\tau_\beta(i_k)}$ where $(i_j)$ are the
nondecreasing sequences, is $h_k(X_n)$.
\qed


A slightly more general result can be derived from the previous
considerations: if one considers a noncommutative ordered alphabet $A_n =
\{a_0<\ldots<a_n\}$ instead of $X_n$, so that one associates with a
sub-diagonal word $c$ the word
$a_c = a_{c_1} a_{c_2} \cdots a_{c_n}$,
the noncommutative series generalizing Equation~(\ref{scode-X}) reads
\begin{equation}
\label{scode-A}
\sum_{\sigma\ \in\ \id_k\! \ssh\beta} a_{\Scode(\sigma)}
= S_k(A'_n) a_{\Scode(\beta)},
\end{equation}
where $A'_n$ is the ordered alphabet on $\{a_0,\ldots,a_n\}$ where $a_i<a_j$
if $i$ is to the left of $j$ in $\tau_S(\beta)$. This property directly comes
from Equation~(\ref{egalite}).
%

\section{Noncommutative and commutative generating function for codes}

We have already mentioned that the evaluations of the scodes and of the
invcodes are the same over the symmetric group, which is obvious since both
are codes and hence run over the set of sub-diagonal words.
This proves the observation of Lascoux.
Actually, the scode and the invcode have much more in common.
The key result is that the sorted vectors $\Scode(\sigma)^\uparrow$ and
$\Invcode(\sigma)^\uparrow$ have the same distribution on inverse descent
classes, a property also shared by the majcode (see
Section~\ref{sec-majcode}).
Equation~(\ref{sg-Sc}) of Theorem~\ref{scode-thsg} gives the closed expression
of the generating function of these statistics.

\begin{theorem}
\label{invcode-thsg}
Let $I=(i_1,\ldots,i_r)$ be a composition.
The sum of the invcodes of the inverses of the elements of $D_{\leq I}$ are
given by the noncommutative generating function
\begin{equation}
\begin{split}
\label{sg-Ic}
F_S(I) := \sum_{\sigma\in D_{\leq I}} a_{\Invcode(\sigma^{-1})} &=
\sum_{\sigma\ \in\ \id_{i_1}\!\!\ssh\cdots\ssh\,\id_{i_r}} a_{\Invcode(\sigma)}
\\ &= S_{i_1}(A_{i_2+\cdots+i_r})\
S_{i_2}(A_{i_3+\cdots+i_r}) \cdots S_{i_{r-1}}(A_{i_r})\ S_{i_r}(A_0).
\end{split}
\end{equation}
This right-hand side will be denoted by $S^I(\A_I)$, where $\A_I$ denotes the
\emph{flag of alphabets}
\begin{equation}
(A_{n-i_1},A_{n-i_1-i_2},\ldots,A_{n-n}).
\end{equation}
\end{theorem}

\Proof
The proof proceeds by induction on the number of parts of $I$. If $I$ has one
part, the only permutation is the identity and the statement is obvious.
Assuming the result for the composition $(i_2,\ldots,i_r)$, let us prove it
for $I$.
First, let $\sigma$ be an element of $\id_{i_2}\!\ssh\cdots\ssh\id_{i_r}$
and let $\gamma$ be any element in $\id_{i_1}\!\ssh\,\sigma$.
Then, $\Invcode_{i_1+k}(\gamma)=\Invcode_k(\sigma)$ for all $k$.
Moreover, the sequence $\Invcode_{k}(\gamma)$ for $k\in[1,i_1]$ is
nondecreasing, since $1,\ldots,i_1$ are in this order in $\gamma$,
and it is bounded by the number of letters of $\sigma$, that is,
$i_2+\cdots+i_r$.

Since the invcode is a bijection, no two words $\gamma$ can have the same
code, hence the same first $i_1$ values, since the other ones are identical.
Finally, the number of elements in $\id_{i_1}\!\ssh\ \sigma$ is equal to the
number of nondecreasing sequences of size $k$ in $[0,i_2+\cdots+i_r]$
(a binomial coefficient), so that all sequences appear, and the sum of
the invcodes of all elements in $\id_{i_1}\!\ssh\,\sigma$ is
$S_{i_1}(A_{i_2+\dots+ i_r}) a_{\Invcode(\sigma)}$.
\qed


\begin{corollary}
The invcodes of the permutations in an inverse descent class are given by
\begin{equation}
\sum_{\sigma\in D_I} a_{\Invcode(\sigma^{-1})}
= \sum_{J\leq I} (-1)^{l(I)-l(J)} S^J(\A_J) =: R_I(\A_I).
\end{equation}
\end{corollary}

Taking the commutative image $(a_i\to x_i)$, we recover the following
expression (see~\cite{Loth2}, chap. 11):

\begin{corollary}
The commutative generating series for the codes on a descent class is given by
the following determinant (a flagged ribbon Schur function)
\begin{equation}
r_I(\X_I) :=
\left|
\begin{matrix}
h_{i_1}(X_{n-i_1}) & h_{i_1+i_2}(X_{n-i_1-i_2}) & \cdots &
 & h_{i_1+\cdots+i_r}(X_{0}) \\
          1        & h_{i_2}(X_{n-i_1-i_2}) & \cdots &
 & h_{i_2+\cdots+i_r}(X_{0}) \\
                   &            1           & \ddots &
 & \vdots \\
                   &                        & \ddots &
 & \\
                   &                        &        & 1
 & h_{i_r}(X_0) \\
\end{matrix}
\right|
\end{equation}
\end{corollary}

It is interesting to observe that this flagged Schur function is in fact a
Schubert polynomial~\cite{Las,L2}.
For example, with $I=(5,1,2)$, one gets the following determinant:
\begin{equation}
r_{512}(X_3,X_2,X_0) =
\left|
\begin{matrix}
h_5(X_3) & h_6(X_2) & h_8(X_0)\\
1        & h_1(X_2) & h_3(X_0)\\
0        & 1        & h_2(X_0)\\
\end{matrix}
\right|
= Y_{20150000}.
\end{equation}

\begin{corollary}
\label{inv-cor}
The commutative generating series of $\Scode(\sigma^{-1})$ on a descent class
coincides with that of $\Invcode(\sigma^{-1})$:
\begin{equation}
\sum_{\sigma\in D_I} x_{\Scode(\sigma^{-1})} =
\sum_{\sigma\in D_I} x_{\Invcode(\sigma^{-1})}
= r_I(\X_I).
\end{equation}
In other words, the sorted vectors $\Invcode(\sigma^{-1})^\uparrow$ and
$\Scode(\sigma^{-1})^\uparrow$ have the same distribution on a descent class.
\end{corollary}

For example, (\ref{2112}),~(\ref{2112-cod}),~(\ref{2112-scod})
present the 19 permutations with descent composition $(2,1,1,2)$ and their
the invcodes and inverse scodes. Both statistics give the sequences of
(\ref{2112-sort}) when sorted.
{\footnotesize
\begin{equation}
\label{2112}
\begin{split}
  & 154326, 164325, 165324, 165423, 254316, 264315, 265314, 265413, 354216,
  364215, \\
  & 365214, 365412, 453216, 463215, 465213, 465312, 563214, 564213, 564312.
\end{split}
\end{equation}

\begin{equation}
\label{2112-cod}
\begin{split}
  & 032100, 042100, 043100, 043200, 132100, 142100, 143100, 143200, 232100,
  242100, \\
  & 243100, 243200, 332100, 342100, 343100, 343200, 442100, 443100, 443200.
\end{split}
\end{equation}

\begin{equation}
\label{2112-scod}
\begin{split}
  & 043200, 013200, 041200, 043100, 243200, 213200, 241200, 243100, 343200,
  313200, \\
  & 341200, 143100, 443200, 413200, 141200, 343100, 113200, 441200 443100.
\end{split}
\end{equation}

\begin{equation}
\label{2112-sort}
\begin{split}
  & 000123, 000124, 001123, 000134, 001124, 001223, 000234, 001134, 001224,
  001233, \\
  & 001234, 001234, 001234, 001244, 001334, 002234, 001344, 002334, 002344.
\end{split}
\end{equation}}%
Note that one can reinforce the parallel between the scode and the invcode by
observing that the sequence of lemmas used in the scode case could be directly
translated in the invcode one, the only modification being the definition of
$\tau_S$, which would be here
\begin{equation}
\label{tau-inv}
\tau_I(\beta) = (0,1,\ldots,n).
\end{equation}

The fact that $\tau_I$ does not depend on $\beta$ and satisfies the same
property as in Equation~(\ref{egalite}) is a satisfactory way to explain why
the invcode has a noncommutative formula over the shuffle classes whereas
$\tau_S$ has not. We shall see in the next Section that
Equation~(\ref{egalite}) does not hold for the $\tau_M$, the permutation
corresponding to the majcode and that this property explains why there is no
noncommutative formula for this code.

\section{Generating function for majcodes}
\label{sec-majcode}

\begin{theorem}
\label{maj-thsg}
The commutative generating function of $\Majcode(\sigma^{-1})$ on a descent
class coincide with that of $\Invcode(\sigma^{-1})$:
\begin{equation}
\sum_{\sigma\in D_I} x_{\Majcode(\sigma^{-1})} = \sum_{\sigma\in D_I}
x_{\Invcode(\sigma^{-1})}
= r_I(X_I).
\end{equation}
In other words, the sorted vectors $\Majcode(\sigma^{-1})^\uparrow$ and
$\Invcode(\sigma^{-1})^\uparrow$ have the same distribution on a descent
class.
\end{theorem}

For example,~(\ref{2112-maj}) presents the majcodes of the 19 permutations
of~(\ref{2112}).
This statistic also gives the sequences of (\ref{2112-sort}) when sorted.
{\footnotesize
\begin{equation}
\label{2112-maj}
\begin{split}
  & 332100, 342100, 343100, 343200, 232100, 242100, 243100, 443200, 132100,
  142100, \\
  & 443100, 243200, 032100, 442100, 143100, 143200, 042100, 043100, 043200.
\end{split}
\end{equation}
}
By inclusion-exclusion, the statement of the theorem is equivalent to
\begin{equation}
\sum_{\sigma\ \in\ \id_{i_1}\!\!\ssh\cdots\ssh\, \id_{i_r}} x_{\Majcode(\sigma)} =
\sum_{\sigma\ \in\ \id_{i_1}\!\!\ssh\cdots\ssh\, \id_{i_r}} x_{\Invcode(\sigma)}
\end{equation}
for all compositions $I=(i_1,\ldots,i_r)$.
This result is a consequence of the next four lemmas.

\begin{lemma}
Let $\beta\in\SG_n$. Then
\begin{equation}
x_{\Majcode(1\ssh_i \beta)} = x_{\tau_M(\beta)(i)}\ x_{\Majcode(\beta)},
\end{equation}
where $\tau_M(\beta)$ is the permutation of $\{0,\ldots,n\}$ defined by

\begin{equation}
\label{tau-maj}
\begin{split}
\tau_M(\beta)(i)=\des(\beta) - j&\quad
   \text{if $i$ is the $j$-th descent of $\beta$,} \\
\tau_M(\beta)(i)=\des(\beta) + j-1&\quad
   \text{if $i$ is the $j$-th rise of $\beta$},
\end{split}
\end{equation}
where the first position is considered as a rise, and $\des$ denotes the
number of descents.
\end{lemma}
For example,
\begin{equation}
\begin{split}
\beta &=\,\, 9\ 4\ 1\ 6\ 2\ 5\ 7\ 3\ 8 \\
\tau_M(\beta) &= 4\ 3\ 2\ 5\ 1\ 6\ 7\ 0\ 8\ 9
\end{split}
\end{equation}
Note that $\tau_M(\beta)$ only depends on the descent composition $J$ of
$\beta$ so that we can also denote it by $\tau_M(J)$.

\Proof
Let $\beta'=1\ssh_i\beta$. It is obvious that $m_{i+1}(\beta')=m_i(\beta)$ for
$i\in [1,n]$. So $x_{\Majcode(\beta')}/x_{\Majcode(\beta)} = x_{m_1(\beta')}$.
Now, if $i$ is a descent of $\beta$, the insertion of $1$ shifts by one
position all the descents to its right.
So, in this case, $x_{m_1(\beta')} = x_{\tau_M(\beta)(i)}$.
Otherwise, if $i=0$, it shifts all descents of $\beta$ by one, so the same
formula holds. Finally, between a rise $i$ and the next rise $j$ in $\beta$,
we have $x_{m_1(1\ssh_j\beta)}-x_{m_1(1\ssh_i\beta)}=1$ since the first one
creates a descent at position $j$ instead of position $i$ and shifts $j-i-1$
descents ($j$ is the \emph{next} rise after $i$) less than $1\ssh_i\beta$.
So the formula holds for all $i$.
\qed

\begin{lemma}
Let $\beta\in\SG_n$ and let $k$ be an element of $[0,n]$.
The set
$\{\tau_M(\beta)(i)\}_{i\in[0,k]}$
is the interval $[d,d+k]$, where $d$ is the number of descents of $\beta$
greater than $k$.
\end{lemma}

\Proof
By construction of $\tau_M(\beta)$, it is obvious that this set is an interval
of $k+1$ elements. Since its smallest element is the number of descents of
$\beta$ greater than $k$, obtained either at the last descent of $\beta$
smaller than or equal to $k$ or at the first rise of $\beta$, we are done.
\qed

From this, we get immediately:

\begin{lemma}
\label{maj-l3}
Let $\beta\in\SG_n$, let $k$ be an element of $[0,n]$
and let $\beta'=1\ssh_k\beta$.
Then,
\begin{equation}
\label{maj-ssens}
\{\tau_M(\beta')(i)\ |\ i\in[0,k]\} = \{\tau_M(\beta)(i)\ |\ i\in[0,k]\}.
\end{equation}
\end{lemma}
The following tableau presents a permutation $\beta$ and its permutation
$\tau_M(\beta)$ and then $1\ssh_5\beta$ and $1\ssh_6\beta$ and their
results by application of $\tau_M$.
\begin{equation}
\label{maj-ex}
\begin{array}{|l|l|l|}
\hline
             & Permutation & \tau_M     \\
\hline
\beta        & 72451836    & 324516078  \\
\hline
1\ssh_5\beta & 835621947   & 4356217089 \\
\hline
1\ssh_6\beta & 835629147   & 3245160789 \\
\hline
\end{array}
\end{equation}


\begin{lemma}
\label{maj-l4}
Let $\beta\in\SG_n$ and $k$ be an integer.
The set of the sorted $k$ first components of the majcodes of the elements in
$id_k \ssh\beta$ is the set of all sequences
$(0\leq i_1\leq i_2\leq\cdots\leq i_k\leq n)$.
In particular, we have
\begin{equation}
\sum_{\sigma\ \in\ \id_k\!\ssh\beta} x_{\Majcode(\sigma)}
 = h_k(X_n) \ \ x_{\Majcode(\beta)}.
\end{equation}
\end{lemma}

\Proof
The first part in proved by induction on $k$. For $k=1$, the result is
obvious. Assume that it is true for $k-1$. We prove that for all
nondecreasing sequences $(0\leq i_1\leq i_2\leq\cdots\leq i_k\leq n)$,
exactly one of its permutations can be obtained as the $k$ first components of
the majcode of an element in $id_k\ssh\beta$. Consider a permutation satisfying
this property. Then the $k$-th component of its majcode is necessarily the
value $i_j$ that is in the rightmost position in $\tau_M(\beta)$ since
otherwise this value could not be anywhere else in the majcode of this
permutation thanks to Equation~(\ref{maj-ssens}). It is indeed the case
of $\sigma=175239486\in 1234\ssh 31542$ whose majcode is
$424133110$, whose fourth component is $1$, which is indeed to the right
of $2$ and $4$ in $\tau_M(31542)=324105$.

We can then conclude by induction, since there are as many nondecreasing
sequences of size $k$ from $0$ to $n$ as elements in the shuffle
$id_k\ssh\beta$.
Then, since $\tau_G$ is a permutation of $[0,n]$, the commutative
image of the sum of the $k$ first components of the majcodes is
$h_k(X_n)$.
\qed

In the case of the majcode, no noncommutative formula holds, even for
$id_k \ssh\beta$. Indeed, there is no fixed order on the set
$\{0,\ldots,n\}$ as in the case of the scode associated with $\beta$. For
example, $1423$ belongs to $123\ssh1$ and its scode is $101\,0$ so that $1<0$
and $0<1$, which is impossible.

\section{Codes and Euler-Mahonian statistics}

Let us say that a general code $\Gcode$ is \emph{compatible with the shuffle}
if the codes of $1\ssh\beta$ are of the form $(i,\Gcode(\beta))$. Since the
code is a bijection, one can then define $\tau(\beta)$ as the permutation of
$\{0,\ldots,n\}$ obtained by sending $i$ to the first value of the code of
$1\ssh_i\beta$.
For example, the scode, the invcode, and the majcode are compatible with the
shuffle, the corresponding permutations $\tau(\beta)$ being respectively
as defined in Equations~(\ref{tau-sco}),~(\ref{tau-inv}),~(\ref{tau-maj}).
We then say that a general code $\Gcode$ compatible with the shuffle is an
\emph{acceptable code} if, for all permutations $\beta$ and all
integers $k$,
\begin{equation}
\label{tau-gen}
\{\tau_G(\beta')(i)\ |\ i\in[0,k]\} = \{\tau_G(\beta)(i)\ |\ i\in[0,k]\}.
\end{equation}
where $\beta'=1\ssh_k\beta$.
For example, the scode, the invcode, and the majcode are acceptable codes, as
it is trivial for the invcode and already done in Lemma~\ref{sco-l3} for the
scode and in Lemma~\ref{maj-l3} for the majcode.

The following Lemma is proved exactly as Lemma~\ref{maj-l4}, by induction on
$k$.

\begin{lemma}
\label{gen-l4}
Let us consider an acceptable code $\Gcode$.

Let $\beta\in\SG_n$ and $k$ be an integer.
Then consider the set of the sorted $k$ first components of the Gcodes of the
elements in $id_k \ssh\beta$.
This set is exactly the set of all sequences
\begin{equation}
(0\leq i_1\leq i_2\leq\cdots\leq i_k\leq n).
\end{equation}

In particular, we have
\begin{equation}
\sum_{\sigma\ \in\ \id_k\! \ssh\, \beta} x_{\Gcode(\sigma)}
= h_k(X_n)\ x_{\Gcode(\beta)}.
\end{equation}
\end{lemma}
This implies the following theorem, which contains Theorems~\ref{scode-thsg}
and~\ref{inv-cor}, and Corollary~\ref{maj-thsg} as particular cases.

\begin{theorem}
\label{gen-tau}
Let us consider an acceptable code $\Gcode$. Then
the commutative generating series of $\Gcode(\sigma^{-1})$ on a descent
class coincides with that of $\Invcode(\sigma^{-1})$:
\begin{equation}
\sum_{\sigma\in D_I} x_{\Gcode(\sigma^{-1})} = \sum_{\sigma\in D_I}
x_{\Invcode(\sigma^{-1})}
= r_I(X_I).
\end{equation}
In other words, the sorted vectors $\Gcode(\sigma^{-1})^\uparrow$ and
$\Invcode(\sigma^{-1})^\uparrow$ have the same distribution on a descent
class.
\end{theorem}

\begin{corollary}
Let $\Gcode$ be an acceptable code.
Then the bi-statistic
\begin{equation}
\left(\sum_{i\in \Gcode(\sigma^{-1})}{i},\qquad \des(\sigma)\right)
\end{equation}
is Euler-Mahonian.
\end{corollary}

\section{An equivalence related to sorted codes}

The previous sections showed the importance of the sorted codes.
We present here a simple construction relating all permutations having a given
sorted Lcode.
Let us say that two words $u$ and $v$ are \emph{L-adjacent} if there exists
four words $w_1,w_2,w_3,w_4$ and three letters $a<b<c$ such that
\begin{equation}
\begin{split}
u =& w_1\, a\, w_2\, c\, w_3\, b\, w_4, \\
v =& w_1\, b\, w_2\, a\, w_3\, c\, w_4, \\
\end{split}
\end{equation}
where all letters of $w_2$ are greater than $b$, and all letters of $w_3$ and
$w_4$ are either smaller than $b$ or greater than $c$.
The \emph{L-equivalence} is the transitive closure of the relation of
L-adjacency. That is, two words $u,v$ are L-equivalent if there exists a
chain of words
$u=u_1, u_2, \dots, u_k=v$,
such that $u_i$ and $u_{i+1}$ are L-adjacent for all $i$. In this case, we
write $u \sim v$.
For example, the L class of $w=31452$ is the set
$\{13542$, $14352$, $21543$, $23514$, $24153$, $24315$, $31452$, $32154$,
$ 32415\}$.
This is \emph{not} a congruence on $A^*$.

\begin{lemma}
\label{l1}
If $u$ and $v$ are L-adjacent, the sorted Lcode of $u$ and $v$ are equal.
\end{lemma}

\Proof
It is easy to check that the parts of the code corresponding to the subwords
$w_i$ are the same in the code of $u$ and in the code of $v$. Moreover, the
part of the code corresponding to $a$ (resp. $b$, $c$) in $u$ is the same as
the part of the code corresponding to $a$ (resp. $c$, $b$) in $v$.
\qed

\noindent
For example $738694152$ and $758634192$ are L-adjacent ($a=3$, $b=5$,
$c=9$) and their Lcodes are
$\Lcode(u) = 625442010 \text{\ and\ } \Lcode(v) = 645422010$.

\begin{lemma}
\label{l2}
Let $w$ be a word avoiding the pattern $132$. Then there exists $w'$ greater
than $w$ in the lexicographic order such that $w$ and $w'$ are L-adjacent.
\end{lemma}

\Proof
Let $w$ be a word containing the pattern $132$. Then there exists an
occurrence of the pattern $acb$ such that the difference $c-b$ is minimal.
Then take $a$ as the rightmost element smaller than $b$ to the left of $c$.
This element is L-adjacent to the element obtained by changing $a$ into $b$,
$c$ into $a$ and $b$ into $c$.  Indeed, by construction, all the words $w_2$,
$w_3$, and $w_4$ satisfy the hypotheses of L-adjacency.
\qed

\begin{proposition}
Two words $u$ and $v$ are L-equivalent iff $u$ and $v$ have same sorted
Lcode.
\end{proposition}

\Proof
The number of sorted Lcodes is a Catalan number, also counting
the permutations avoiding $132$. By Lemma~\ref{l1}, the number of L classes is
greater than or equal to the Catalan number. By Lemma~\ref{l2}, the number of
L classes is smaller than or equal to the same Catalan number. So the number
of L classes is given by the same Catalan number and the classes are the same
as the classes of words by sorted Lcode.
\qed

Hence, each L class contains exactly one permutation avoiding $132$ (its
greatest element). In the same way, one can prove that each L class contains
exactly one permutation avoiding $213$ (its smallest element). The map sending
a permutation to the maximal element of its L-class consists in sorting the
Lcode of the permutation, whereas the map sending a permutation to the minimal
element of its L-class consists in building, from right to left, the
permutation of its Lcode that has the greatest possible value in each
position.
For example, with $\sigma = 682547193$ whose Lcode is $561322010$, the maximal
element of its L-class has Lcode $653221100$, so is permutation $764352819$,
and the minimal element of its L-class has Lcode $016523210$, so is
$139857642$.

\section{Tables}

We first give the values of the three statistics used in the paper, namely,
the inverse Lehmer code, the majcode and the scode.
In the table, the permutations are sorted by inverse descent classes.

\begin{equation*}
\begin{array}{cccc@{\qquad\qquad}cccc}
  \sigma & \Invcode & \Majcode & \Scode    & \sigma & \Invcode & \Majcode
 & \Scode  \\
   1234 & 0000 & 0000 & 0000   &   4321 & 3210 & 3210 & 3210\\
\\
   1243 & 0010 & 1110 & 0010   &   3214 & 2100 & 2100 & 3200\\
   1423 & 0110 & 1010 & 0110   &   3241 & 3100 & 3100 & 1200\\
   4123 & 1110 & 0010 & 1110   &   3421 & 3200 & 3200 & 3100\\
\\
   1324 & 0100 & 1100 & 0200   &   2143 & 1010 & 2110 & 3010\\
   1342 & 0200 & 1200 & 0100   &   2413 & 2010 & 0110 & 1010\\
   3124 & 1100 & 0100 & 2200   &   2431 & 3010 & 3110 & 2010\\
   3142 & 1200 & 2200 & 2100   &   4213 & 2110 & 2010 & 3110\\
   3412 & 2200 & 0200 & 1100   &   4231 & 3110 & 3010 & 2110\\
\\
   1432 & 0210 & 2210 & 0210   &   2134 & 1000 & 1000 & 3000\\
   4132 & 1210 & 1210 & 1210   &   2314 & 2000 & 2000 & 2000\\
   4312 & 2210 & 0210 & 2210   &   2341 & 3000 & 3000 & 1000\\
\end{array}
\end{equation*}


The next expressions present the first values of polynomials $r_I(\X_I)$ where
the monomial $x_0^{i_0}\cdots x_k^{i_k}$ is represented by
$[0^{i_0},1^{i_1},2^{i_2},\ldots,k^{i_k}]$.

{
\footnotesize
\begin{alignat*}{1}
r_{3}   & =  [000] \\
r_{21}   & =  [001] + [011] \\
r_{12}   & =  [001] + [002] \\
r_{111}   & =  [012] \\
\end{alignat*}

\begin{alignat*}{1}
r_{4}   & =  [0000] \\
r_{31}   & =  [0001] + [0011] + [0111] \\
r_{22}   & =  [0001] + [0002] + [0011] + [0012] + [0022] \\
r_{211}   & =  [0012] + [0112] + [0122] \\
r_{13}   & =  [0001] + [0002] + [0003] \\
r_{121}   & =  [0011] + [0012] + [0013] + [0112] + [0113] \\
r_{112}   & =  [0012] + [0013] + [0023] \\
r_{1111}   & =  [0123] \\
\end{alignat*}

\begin{alignat*}{1}
r_{5}   & =  [00000] \\
r_{41}   & =  [00001] + [00011] + [00111] + [01111] \\
r_{32}   & =  [00001] + [00002] + [00011] + [00012] + [00022] + [00111] + [00112] + [00122] + [00222] \\
r_{311}   & =  [00012] + [00112] + [00122] + [01112] + [01122] + [01222] \\
r_{23}   & =  [00001] + [00002] + [00003] + [00011] + [00012] + [00013] + [00022] + [00023] + [00033] \\
r_{221}   & =  [00011] + [00012] + [00013] + [00111] + 2\ [00112] + 2\ [00113] + [00122] + [00123] + [00133] \\
          &\ + [01112] + [01113] + [01122] + [01123] + [01133] \\
r_{212}   & =  [00012] + [00013] + [00023] + [00112] + [00113] + [00122] + 2\ [00123] + [00133] \\
          &\ + [00223] + [00233] \\
r_{2111}   & =  [00123] + [01123] + [01223] + [01233] \\
r_{14}   & =  [00001] + [00002] + [00003] + [00004] \\
r_{131}   & =  [00011] + [00012] + [00013] + [00014] + [00111] + [00112] + [00113] + [00114] + [01112] \\
          &\ + [01113] + [01114] \\
r_{122}   & =  [00011] + 2\ [00012] + [00013] + [00014] + [00022] + [00023] + [00024] + [00112] + [00113] \\
          &\ + [00114] + [00122] + [00123] + [00124] + [00223] + [00224] \\
r_{1211}   & =  [00112] + [00122] + [00123] + [00124] + [01122] + [01123] + [01124] + [01223] + [01224] \\
r_{113}   & =  [00012] + [00013] + [00014] + [00023] + [00024] + [00034] \\
r_{1121}   & =  [00112] + [00113] + [00114] + [00123] + [00124] + [00134] + [01123] + [01124] + [01134] \\
r_{1112}   & =  [00123] + [00124] + [00134] + [00234] \\
r_{11111}   & =  [01234] \\
\end{alignat*}

}


\end{document}